\documentclass[12pt,reqno]{amsart}

\usepackage[latin1]{inputenc}
\usepackage[T1]{fontenc}
\usepackage{amssymb,rcs,url}

\theoremstyle{plain}
\newtheorem{theorem}{Theorem}[section]
\newtheorem{lemma}[theorem]{Lemma}
\newtheorem{proposition}[theorem]{Proposition}

\theoremstyle{definition}

\theoremstyle{remark}
\newtheorem{remark}[theorem]{Remark}
\newtheorem*{acknowledgements}{Acknowledgements}

\numberwithin{equation}{section}

\newcommand{\Lie}[1]{\textsl{#1}}
\newcommand{\lie}[1]{\mathfrak{#1}}

\DeclareMathOperator{\ad}{ad}

\DeclareMathOperator{\CP}{\mathbb CP}

\DeclareMathOperator{\HP}{\mathbb HP}

\DeclareMathOperator{\Id}{Id}
\DeclareMathOperator{\im}{Im}
\newcommand{\inp}[3][]{\left\langle #2,#3\right\rangle_{#1}}

\DeclareMathOperator{\LSP}{\lie{sp}}

\newcommand{\norm}[2][]{\left\lVert #2\right\rVert_{#1}}

\DeclareMathOperator{\re}{Re}

\DeclareMathOperator{\Sl}{\lie{sl}}
\DeclareMathOperator{\SO}{\Lie{SO}}
\DeclareMathOperator{\so}{\lie{so}}
\DeclareMathOperator{\SP}{\Lie{Sp}}

\newcommand{\Span}[2][\mathbb C]{\operatorname{Span}_{#1}{\{#2\}}}
\DeclareMathOperator{\Spin}{\Lie{Spin}}

\DeclareMathOperator{\SU}{\Lie{SU}}
\DeclareMathOperator{\su}{\lie{su}}

\DeclareMathOperator{\Un}{\Lie U}

\newenvironment{spmatrix}{\left(\smallmatrix}{\endsmallmatrix\right)}

\def\sddots{\mathinner{\mkern1mu\raise1ex
    \vbox{\kern0.2ex\hbox{.}}\mkern2mu
    \raise0.5ex\hbox{.}\mkern2mu\hbox{.}\mkern1mu}}

\makeatletter

\newcommand{\ltwoa}{\mathord{\mathpalette\tw@a<}}
\newcommand{\rtwoa}{\mathord{\mathpalette\tw@a>}}
\newcommand{\tw@a}[2]{\ooalign{\hfil$#1 #2$\hfil\crcr
$#1 \mathrel=\joinrel\mathrel=$\crcr}}

\newcommand{\lthreea}{\mathord{\mathpalette\thr@a<}}
\newcommand{\rthreea}{\mathord{\mathpalette\thr@a>}}
\newcommand{\thr@a}[2]{\ooalign{\hfil$#1 #2$\hfil\crcr
$#1 \mathrel=\joinrel\mathrel=$\crcr
$#1 \relbar\joinrel\relbar$}}

\newcommand{\Dyatop}[2]{\mathord{\mathpalette{\Dy@top{#1}{#2}}{}}}
\newcommand{\Dy@top}[3]{\overset{#3 #1}{#3 #2}}
\makeatother

\renewcommand{\today}{%
  \relax\number\day\space
  \ifcase\month\or
  January\or February\or March\or April\or May\or June\or
  July\or August\or September\or October\or November\or December\fi
  , \number\year}

\RCS $Revision: 1.9 $
\RCS $Date: 1999/06/23 12:03:26 $

\begin{document}

\title{The HyperKähler Geometry Associated to~Wolf Spaces}

\author{Piotr Kobak}
\address[Kobak]{Insytut matematyki\\
Uniwersytet Jagiello\'nski\\
ul.\ Reymon\-ta 4\\
30-059 Kraków\\
Poland}
\email{kobak@im.uj.edu.pl}

\author{Andrew Swann}
\address[Swann]{Department of Mathematical Sciences\\
University of Bath\\
Claverton Down\\ 
Bath BA2 7AY\\
England}
\email{A.F.Swann@maths.bath.ac.uk}
\curraddr{From 1/9/99: Department of Mathematics and Computer Science,
SDU, Odense University, Campusvej
55, DK-5230 Odense M, Denmark} 
\email{swann@imada.sdu.dk}

\maketitle

\section{Introduction}
\label{sec:introduction}
One of the glories of homogeneous geometry is Cartan's classification of
the compact Riemannian symmetric spaces
\cite{Cartan:symmetric-1,Cartan:symmetric-2}.  Many manifolds that play a
central rôle in geometry are symmetric and it is fascinating to look for
patterns in the presentations~$G/H$.  One obvious family is provided by the
sphere~$S^n=\SO(n+1)/\SO(n)$, complex projective
space~$\CP(n)=\Un(n+1)/(\Un(n)\Un(1))$, quaternionic projective
space~$\HP(n)=\SP(n+1)/(\SP(n)\SP(1))$ and the Cayley projective
plane~$\Lie F_4/\Spin(9)$.  Another consists of the Hermitian symmetric
spaces: these are of the form $G/(\Un(1)L)$ (see
\cite{Burstall-Rawnsley:twistor}).  However, the most surprising is the
family of quaternionic symmetric spaces $W(G):=G/(\SP(1)K)$, which has the
feature that there is precisely one example for each compact simple
simply-connected Lie group~$G$.  The manifolds in this last family have
become known as Wolf spaces following~\cite{Wolf:quaternionic}.
Alekseevsky \cite{Alekseevsky:compact} proved that they are the only
homogeneous positive quaternionic Kähler manifolds
(cf.~\cite{Alekseevsky-C:unimodular}).

Wolf showed that the quaternionic symmetric spaces may be constructed by
choosing a highest root~$\alpha$ for~$\lie g^{\mathbb C}$.  The
corresponding root vector~$E_\alpha$ is a nilpotent element in~$\lie
g^{\mathbb C}$.  In \cite{Swann:MathAnn} it was shown that there is a
fibration of the nilpotent adjoint orbit $\mathcal
O_{\text{min}}=G^{\mathbb C}\cdot E_\alpha$ over the Wolf space~$W(G)$.

Nilpotent orbits~$\mathcal O$ in~$\lie g^{\mathbb C}$ have a rich and
interesting geometry.  Firstly, they are complex submanifolds of~$\lie
g^{\mathbb C}$ with respect to the natural complex structure~$I$.
Secondly, the construction of Kirillov, Kostant and Souriau endows them
with a $G^{\mathbb C}$-invariant complex symplectic form~$\omega_c$.  It is
natural to ask whether one can find a metric making the orbit hyperKähler,
i.e., can one find a Riemannian metric~$g$ on~$\mathcal O$, such that the
real and imaginary parts of~$\omega_c$ are Kähler forms with respect to
complex structures $J$ and~$K$ satisfying $IJ=K$.  By identifying $\mathcal
O$ with a moduli space of solutions to Nahm's equations,
Kronheimer~\cite{Kronheimer:nilpotent} showed that there is indeed such a
hyperKähler metric on~$\mathcal O$.  This hyperKähler structure is
invariant under the compact group~$G$, and has the important additional
property that it admits~\cite{Swann:MathAnn} a hyperKähler
potential~$\rho$: a function that is simultaneously a Kähler potential with
respect to $I$, $J$ and~$K$.  Using~$\rho$, one can define an action
of~$\mathbb H^*$ on~$\mathcal O$ such that the quotient is a quaternionic
Kähler manifold.  It is in this way that one may obtain the Wolf
space~$W(G)$ from~$\mathcal O_{\text{min}}$.  In contrast to the
semi-simple case \cite{Bielawski:homogeneous}, currently one does not know
how many invariant hyperKähler metrics a given nilpotent orbit admits.

The aim of this paper is to study the hyperKähler geometry of~$\mathcal
O_{\text{min}}$ in an elementary way.  We look for all hyperKähler metrics
on $\mathcal O_{\text{min}}$ with a $G$-invariant Kähler potential and
which are compatible with the complex symplectic structure.  Note that we
do not restrict our attention to metrics with hyperKähler potentials.  We
derive a simple formula for the a priori unknown complex structure~$J$.
The orbit $\mathcal O_{\text{min}}$ is particularly straight-forward to
study in this way, since $G$~acts with orbits of codimension one.  This
means that the metrics we obtain are already known, they are covered by the
classification \cite{Dancer-Swann:hK-cohom1}, but it is interesting to see
how these metrics can be constructed directly form their potentials.  In
agreement with the classification, the hyperKähler structure is found to be
unique, unless $\lie g=\su(2)$, in which case one obtains a one-dimensional
family of metrics, the Eguchi-Hanson metrics.

\begin{acknowledgements}
  We are grateful for financial support from the \textsc{Epsrc} of Great
  Britain and \textsc{Kbn} in Poland.
\end{acknowledgements}

\section{Definitions}
\label{sec:definitions}
On the simple complex Lie algebra~$\mathfrak g^{\mathbb C}$, let
$\inp\cdot\cdot$ be the \emph{negative} of the Killing form and let
$\sigma$~be a real structure giving a compact real form~$\mathfrak g$
of~$\mathfrak g^{\mathbb C}$.  An element~$X$ of~$\lie g^{\mathbb C}$ is
said to be nilpotent if $(\ad_X)^k=0$ for some integer~$k$.  Let $\mathcal
O$ be the orbit of a nilpotent element~$X$ under the adjoint action
of~$G^{\mathbb C}$.  At $X\in\mathcal O$, the vector field generated by~$A$
in~$\mathfrak g^{\mathbb C}$ is $\xi_A=[A,X]$.  Using the Jacobi identity
it is easy to see that these vector fields satisfy $[\xi_A,\xi_B] =
\xi_{-[A,B]}$, for $A,B\in\mathfrak g^{\mathbb C}$.  The orbit $\mathcal
O$~is a complex submanifold of the complex vector space~$\lie g^{\mathbb
C}$ and so has a complex structure~$I$ given by $I\xi_A = i \xi_A =
\xi_{iA}$.

On a hyperKähler manifold~$M$ with complex structures $I$, $J$ and~$K$ and
metric~$g$, we define Kähler two-forms by $ \omega_I(X,Y) = g(X,IY)$, etc.,
for tangent vectors $X$ and~$Y$.  The condition that a function~$\rho\colon
M\to \mathbb R$ be a Kähler potential for~$I$ is
\begin{equation}
  \label{eq:dId}
  \omega_I
  = - i \partial_I \overline{\partial_I} \rho
  = - i d \overline{\partial_I} \rho
  = - \tfrac i2 d (d-iId) \rho
  = - \tfrac 12 d I d \rho.
\end{equation}
On the orbit~$\mathcal O$, the complex symplectic form of Kirillov, Kostant
and Souriau is given by $ \omega_c(\xi_A,\xi_B)_X = \inp X{[A,B]} = -
\inp{\xi_A}B $.  

We will be looking for hyperKähler structures with Kähler potential~$\rho$
and such that $\omega_c=\omega_J +i\omega_K$.  This will be done by
computing the Riemann metric~$g$ defined by~$\rho$ via~\eqref{eq:dId} and
then using this to determine an endomorphism~$J$ of~$T_X\mathcal O$ via
$\omega_J=g(\cdot,J\cdot)$.  The constraints on~$\rho$ will come from the
two conditions that $g$~is positive definite and that $J^2=-1$.

\section{Highest Roots and Minimal Orbits}
\label{sec:minimal}
Choose a Cartan subalgebra~$\lie h$ of~$\lie g^{\mathbb C}$.  Fix a system
of roots~$\Delta$ with positive roots $\Delta_+$.  We write $\lie g_\beta$
for the root space of~$\beta\in\Delta$.  Choose a Cartan basis
$\{E_\beta,H_\beta,F_\beta:\beta\in \Delta_+\}$, which we may assume is
compatible with the real structure~$\sigma$, in the sense that
$\sigma(E_\beta)=-F_\beta$ and $\sigma(H_\beta)=-H_\beta$.  One important
property of the Cartan basis is that for each~$\beta$,
$\Span{E_\beta,H_\beta,F_\beta}$ is a subalgebra of~$\lie g^{\mathbb C}$
isomorphic to~$\Sl(2,\mathbb C)$.

The Lie algebra~$\Sl(2,\mathbb C)$ has Cartan basis
\begin{equation}
  \label{eq:Cartan}
  E= \begin{spmatrix} 0&1\\ 0&0 \end{spmatrix},
  \quad H= \begin{spmatrix} 1&0\\ 0&-1 \end{spmatrix},
  \quad F= \begin{spmatrix} 0&0\\ 1&0 \end{spmatrix}. 
\end{equation}
The irreducible representations of $\Sl(2,\mathbb C)$ are the symmetric
powers $S^k=S^k\mathbb C^2$ of the fundamental representation~$S^1=\mathbb
C^2$.  The representation~$S^k$ has dimension $k+1$ and $E$, $H$ and $F$
act as
\begin{gather}
  \label{eq:Sk}
  \varphi_E=
  \begin{spmatrix}
    0 & 1 & \\
    & 0 & 2 & \\
    & & \sddots & \sddots & \\
    & & & 0 & k  \\
    & & & & 0 
  \end{spmatrix},\quad
  \varphi_H=
  \begin{spmatrix}
    k & \\
    & k-2 & \\
    & & \sddots &\\
    & & & 2-k \\
    & & & & -k
  \end{spmatrix}\\
  \text{and}\quad
  \varphi_F=
  \begin{spmatrix}
    0 &  \\
    k & 0 &  \\
    & \sddots & \sddots & \\
    & & 2 & 0 &   \\
    & & & 1 & 0 
  \end{spmatrix}\notag
\end{gather}
respectively.  In particular, $(\varphi_E)^{k+1}=0$ and $(\varphi_E)^k$ has
rank one, with image the $k$-eigenspace of~$\varphi_H$.

Let $\alpha\in\Delta_+$~be a highest root; this is characterised by the
condition $[E_\alpha,E_\beta]=0$ for all $\beta\in\Delta_+$.  We define
$\mathcal O_{\text{min}}$~to be the adjoint orbit of~$E_\alpha$ under the
action of~$G^{\mathbb C}$.  Define $\Sl(2,\mathbb
C)_\alpha:=\Span{E_\alpha,H_\alpha,F_\alpha}$.  

\begin{proposition}
  \textup{(i)} Under the action of $\Sl(2,\mathbb C)_\alpha$ the Lie
  algebra~$\lie g^{\mathbb C}$ decomposes as
  \begin{equation*}
    \lie g^{\mathbb C} \cong \Sl(2,\mathbb C)_\alpha \oplus \lie k^{\mathbb C}
    \oplus (V\otimes S^1),
  \end{equation*}
  where $\lie k^{\mathbb C}$ is the centraliser of~$\Sl(2,\mathbb C)$,
  $V$~is a $\lie k^{\mathbb C}$-module.
  
  \textup{(ii)} The action of the compact group~$G$ on the nilpotent
  orbit~$\mathcal O_{\text{min}}$ has cohomogeneity one.
\end{proposition}

\begin{proof}
  (i)~Consider the action of~$\ad E_\alpha$ on~$\lie g^{\mathbb C}$.  For
  $\beta\in\Delta_+$, we have $[E_\alpha,F_\beta]\in \lie
  g_{\alpha-\beta}$.  If $\beta\ne\alpha$, then we have two cases: (a)~if
  $\alpha-\beta$ is not a root then $\lie g_{\alpha-\beta}=\{0\}$ and
  $[E_\alpha,F_\beta]=0$; (b)~if $\alpha-\beta$ is a root, then the
  condition that $\alpha$~is a highest root implies
  $\alpha-\beta\in\Delta_+$, since otherwise $\alpha-\beta=-\gamma$ for
  some $\gamma\in\Delta_+$ and then $[E_\alpha,E_\gamma]$~is non-zero,
  which for a highest root~$\alpha$ is impossible.  We therefore have that
  $(\ad E_\alpha)^2$~is zero on the complement of~$\Sl(2,\mathbb C)_\alpha$
  and the decomposition follows.
  
  (ii) At $E_\alpha$ the tangent space to~$\mathcal O_{\text{min}}$ is
  \begin{equation*}
    \ad_{E_\alpha}\lie g^{\mathbb C} = \Span{E_\alpha,H_\alpha} +
    \Span{E_{\alpha-\beta}:\beta\in\Delta_+}.  
  \end{equation*}
  The real Lie algebra~$\lie g$ is the real span of
  $\{E_\beta-F_\beta,iH_\beta,i(E_\beta+F_\beta)\}$.  Thus the tangent
  space~$\ad_{E_\alpha}\lie g$ to the $G$-orbit is 
  \begin{equation*}
    \Span[\mathbb R]{iE_\alpha,H_\alpha,iH_\alpha} +
    \Span[\mathbb R]{E_{\alpha-\beta},iE_{\alpha-\beta} :
    \beta\in\Delta_+} 
  \end{equation*}
  and we see that it has codimension one in~$T_{E_\alpha}\mathcal
  O_{\text{min}}$, the complement being $\mathbb RE_\alpha$.  As $G$~is
  compact, this implies $G$~acts with cohomogeneity one.
\end{proof}

As in \cite{Dancer-Swann:qK-cohom1}, it is possible to use this result to
show that $\mathcal O_{\text{min}}$ is the minimal with respect to the
partial order on nilpotent orbits given by inclusions of closures.  This
explains the name~$\mathcal O_{\text{min}}$, but will not be needed in the
subsequent discussion.

\section{Kähler Potentials in Cohomogeneity One}
Let $\rho\colon\mathcal O_{\text{min}} \to \mathbb R$ be a smooth function
invariant under the action of the compact group~$G$.  The group $G$~acts
with cohomogeneity one, and the function $\eta(X)=\norm X^2 = \inp X{\sigma
X}$ is $G$-invariant and distinguishes orbits of~$G$.  We may therefore
assume that $\rho$~is just a function of~$\eta$, i.e., $\rho = \rho(\eta)$.

We wish to consider $\rho$ as a Kähler potential for the complex
manifold~$(\mathcal O_{\text{min}},I)$.  The corresponding Kähler form is
given by~\eqref{eq:dId}: 
\begin{equation}
  \label{eq:omegaI-min-t}
  \omega_I
  = -\tfrac12d(\rho'I d\eta)
  = -\tfrac12\rho'dId\eta -\tfrac12\rho''d\eta\wedge Id\eta,
\end{equation}
where $\rho'=d\rho/d\eta$, etc.

\begin{lemma}
  The Kähler form defined by $\rho(\eta)$ is
  \begin{equation}
    \label{eq:omegaI-min}
    \omega_I(\xi_A,\xi_B)
    = 2\im         
    \left(
      \rho'\inp{\xi_A}{\sigma\xi_B}
      + \rho'' \inp{\xi_A}{\sigma X}\inp{\sigma \xi_B}X
    \right)
    .
  \end{equation}
\end{lemma}

\begin{proof}
  The exterior derivative of~$\eta$ is
  \begin{equation}
    \label{eq:d-eta}
    d\eta(\xi_A)_X
    = \inp{[A,X]}{\sigma X} + \inp X{\sigma[A,X]}
    = 2 \re \inp{\xi_A}{\sigma X}
  \end{equation}
  so $Id\eta(\xi_A)_X = 2\im \inp{\xi_A}{\sigma X}$ and hence
  \begin{equation*}
    (d\eta\wedge Id\eta)(\xi_A,\xi_B)
    =  - 4 \im \left( \inp{\xi_A}{\sigma X}\inp{\sigma \xi_B}X \right).
  \end{equation*}
  Using the Jacobi identity we find that the exterior derivative
  of~$Id\eta$ is given by
  \begin{equation*}
    \begin{split}
      dId\eta(\xi_A,\xi_B)_X
      &= \xi_A(Id\eta(\xi_B)) - \xi_B(Id\eta(\xi_A)) -
      Id\eta([\xi_A,\xi_B]) \\
      &= 2 \im \inp{\xi_B}{\sigma \xi_A} + 2 \im \inp{[B,[A,X]]}{\sigma X}
      \\
      &\qquad - 2 \im \inp{\xi_A}{\sigma \xi_B} - 2 \im
      \inp{[A,[B,X]]}{\sigma X} \\
      &\qquad + 2 \im \inp{[[A,B],X]}{\sigma X} \\
      &= - 4 \im \inp{\xi_A}{\sigma \xi_B}
    \end{split}
  \end{equation*}
  Putting these expressions into~\eqref{eq:omegaI-min-t} gives the result.
\end{proof}
  
Using the relation $g(\xi_A,\xi_B)=\omega_I(I\xi_A,\xi_B)$, we can now
obtain the induced metric on~$\mathcal O_{\text{min}}$.  In general, this
metric will be indefinite; the signature may be determined by considering
$\Span[\mathbb R]{X,\sigma X}$ and its orthogonal complement with respect
to the Killing form.

\begin{proposition}
  The pseudo-Kähler metric defined by~$\rho(\eta)$ is
  \begin{equation}
    \label{eq:g-min}
    g(\xi_A,\xi_B) = 2\re \left( \rho'\inp{\xi_A}{\sigma\xi_B} +
      \rho'' \inp{\xi_A}{\sigma X}\inp{\sigma \xi_B}X \right) .
  \end{equation}
  This is positive definite if and only if $\rho' >
  \max\{0,-\eta\rho''\}$.\qed 
\end{proposition}

\section{HyperKähler Metrics}
Given a function $\rho(\eta)$ on~$\mathcal O_{\text{min}}$ we have obtained
a metric~$g$.  Let us assume that $g$~is non-degenerate.  Using the
definition of~$\omega_c$ and its splitting into real imaginary parts, we
get endomorphisms $J$ and~$K$ of~$T_X\mathcal O_{\text{min}}$ via
\begin{equation*}
  g(\xi_A,\xi_B) = \omega_J(J\xi_A,\xi_B) = - \re \inp{J\xi_A}B,
\end{equation*}
etc.  This implies that
\begin{equation}
  \label{eq:J-min}
  J_X\xi_A = -2\rho'\left[X,\sigma\xi_A\right]
  -2\rho''\inp{\sigma\xi_A}X\left[X,\sigma X\right] .
\end{equation}
and $K=IJ$.  Note that \eqref{eq:J-min} implies $JI=-K$.  

Suppose $J^2=-1$ and that $g$~is positive definite.  Then we have $I$, $J$
and~$K$ satisfying the quaternion identities, and with $\omega_I$,
$\omega_J$ and~$\omega_K$ closed two-forms.  By a result of
Hitchin~\cite{Hitchin:Montreal}, this implies that $I$, $J$ and~$K$ are
integrable and that $g$ is a hyperKähler metric.

\begin{proposition}
  \label{prop:sl2}
  The nilpotent orbit of~$\Sl(2,\mathbb C)$ has a one-parameter family
  of hyperKähler metrics with $\SU(2)$-invariant Kähler potential and
  compatible with the Kostant-Kirillov-Souriau complex symplectic
  form~$\omega_c$.
\end{proposition}

\begin{proof}
  The algebra $\Sl(2,\mathbb C)$ has only one nilpotent orbit~$\mathcal
  O=\mathcal O_{\text{min}}$ and this has real dimension~$4$.  Using the
  action of~$\SU(2)$ we may assume that $X=tE$, where $t>0$ and $E$~is
  given by~\eqref{eq:Cartan}.  Then $T_X\mathcal O$ is spanned by $H$ and
  $E$.  We have $J_X H = -4\rho't\, E$ and $J_XE = 2t(\rho'+\eta\rho'')H$,
  which implies $J^2=-\Id$ if and only if $8t^2 ({\rho'}^2+\eta\rho'\rho'')
  = 1$.  Now $\eta(E)=4$, so we get the following 
  ordinary differential equation for~$\rho$:
  \begin{equation*}
    2(\eta {\rho'}^2+\eta^2\rho'\rho'')= 1.
  \end{equation*}
  The left-hand side of this equation is $(\eta^2{\rho'}^2)'$, so $
  \rho'=\sqrt{\eta+c}/\eta$, for some real constant~$c$.  For this to be
  defined for all positive~$\eta$, we need $c\geqslant0$.  Now
  $\rho''=-(\eta+2c)/(2\eta^2\sqrt{\eta+c})$, so the metric is 
  \begin{equation}
    \label{eq:g-su2}
    \begin{split}
      g(\xi_A,\xi_B) = \frac 1{\eta^2\sqrt{\eta+c}} \re \bigl(
      &2\eta(\eta+c) \inp{\xi_A}{\sigma\xi_B} \\
      &\quad - (\eta+2c)\inp{\xi_A}{\sigma X}\inp{\sigma \xi_B}X \bigr),
    \end{split}
  \end{equation}
  which is positive definite.
\end{proof}

This hyperKähler metric is of course well-known.  We put it in standard
form as follows.  Using~\eqref{eq:d-eta}, we find
$(\partial/\partial\eta)=E/(8t)$ at $X=tE$.  An $\SU(2)$-invariant basis of
$T_X\mathcal O$ is now given by
$\{\partial/\partial\eta,\xi_{s_1},\xi_{s_2},\xi_{s_3}\}$, where
\begin{equation*}
  s_1=\tfrac12
  \begin{spmatrix}
    0&1\\
    -1&0
  \end{spmatrix},
  \
  s_2=\tfrac12
  \begin{spmatrix}
    0&i\\
    i&0
  \end{spmatrix},
  \
  s_3=\tfrac12
  \begin{spmatrix}
    i&0\\
    0&-i
  \end{spmatrix}
  .
\end{equation*}
This basis is orthogonal with respect to~\eqref{eq:g-su2} and in terms of
the dual basis of one-forms is $\{d\eta,\sigma_1,\sigma_2,\sigma_3\}$,
$g$~is
\begin{equation*}
  \frac1{4\eta^2\rho'}d\eta^2 + \eta\rho'\left(\sigma_1^2+\sigma_2^2\right)
  +\frac1{\rho'}\sigma_3^2.
\end{equation*}
Substituting $\eta=(r/2)^4-c$, we get 
\begin{equation*}
  g=W^{-1}dr^2 + \frac{r^2}4(\sigma_1^2+\sigma_2^2+W\sigma_3^2),
\end{equation*}
with $W=1-16c/r^4$, which are the Eguchi-Hanson
metrics~\cite{Eguchi-H:self-dual}.

\begin{theorem}
  \label{thm:min}
  For $\lie g^{\mathbb C}\ne\Sl(2,\mathbb C)$, the minimal nilpotent
  orbit~$\mathcal O_{\text{min}}$ admits a unique hyperKähler metric with
  $G$-invariant Kähler potential compatible with the complex symplectic
  form~$\omega_c$.
\end{theorem}

\begin{proof}
  Let $\alpha$ be a highest root.  Using the action of~$G$, we may assume
  that $X=tE_\alpha$, for some $t>0$.  On $\xi_A\in\Sl(2,\mathbb
  C)_\alpha$, the condition $J^2=-\Id$ gives $8t^2
  ({\rho'}^2+\eta\rho'\rho'') = 1$, as in Proposition~\ref{prop:sl2}.
  Putting $\lambda^2=\eta(E_\alpha)$, we have $t^2=\eta(X)/\lambda^2$ and
  hence $\rho'=\sqrt{\lambda^2\eta+c}/2\eta$.  Now for $\xi_A$
  Killing-orthogonal to~$\Sl(2,\mathbb C)$, we have
  \begin{equation*}
    J \xi_A =-2\rho'[X,\sigma\xi_A]=-2t\rho'[E_\alpha,\sigma\xi_A]
  \end{equation*}
  and hence
  \begin{equation*}
    J^2 \xi_A = - (4\eta{\rho'}^2/\lambda^2) \ad_{E_\alpha} \ad_{F_\alpha}
    \xi_A = -\left(1+\frac c{\lambda^2\eta} \right) \ad_{E_\alpha}
    \ad_{F_\alpha} \xi_A.  
  \end{equation*}
  As $\eta$~is not constant, the condition $J^2=-\Id$ implies $c=0$ and we
  have a unique hyperKähler metric.
\end{proof}

The proof enables us to write down $J$ explicitly for $\mathcal
O_{\text{min}}$ in $\lie g^{\mathbb C}\ne \Sl(2,\mathbb C)$:
\begin{equation*}
  J_X\xi_A = - \frac\lambda{2\eta^{3/2}}
  \left(
    2\eta[X,\sigma\xi_A] - \inp{\sigma\xi_A}X [X,\sigma X]
  \right).
\end{equation*}
The number $\lambda^2$ is a constant depending only on the Lie
algebra~$\lie g^{\mathbb C}$, with values $2n$ ($\Sl(n,\mathbb C)$,
$\LSP(n-1,\mathbb C)$, $\so(n+2,\mathbb C)$), $8$ ($\Lie G_2$), $18$ ($\Lie
F_4$), $24$ ($\Lie E_6$), $36$ ($\Lie E_7$), $70$ ($\Lie E_8$).

\begin{remark}
  Theorem~\ref{thm:min} only assumes that $\rho$ is a Kähler potential.
  However, the uniqueness result implies that this potential is in fact
  hyperKähler (cf.~\cite{Swann:MathAnn}).  This corresponds to
  Proposition~\ref{prop:sl2}, where $\rho$~is a hyperKähler potential only
  when $c=0$.
\end{remark}

Finally, let us observe that the form of the potential determines the
nilpotent orbit.

\begin{proposition}
  If a nilpotent orbit~$\mathcal O$ has a Kähler potential~$\rho$ that is
  only a function of~$\eta=\norm X^2$ and which defines a hyperKähler
  structure compatible with $\omega_c$, then $\mathcal O$ is a minimal
  nilpotent orbit.
\end{proposition}

\begin{proof}
  Choose $X\in\mathcal O$, such that $\Span{X,\sigma X,[X,\sigma X]}$ is a
  subalgebra isomorphic to~$\Sl(2,\mathbb C)$; this is always possible by a
  result of Borel (cf.~\cite{Kobak-Swann:hk-c2}).  Let $X=tE$, for $t>0$,
  and write $\lie g^{\mathbb C}=\Sl(2,\mathbb C)\oplus\lie m$.  The proofs
  of Proposition~\ref{prop:sl2} and Theorem~\ref{thm:min} imply that
  $\rho'=\lambda\eta^{-1/2}/2$ and $J^2\xi_A=-\ad_E\ad_F\xi_A$ on~$\lie m$.
  Let $S^k$, $k>0$, be an irreducible $\Sl(2,\mathbb C)$-summand of~$\lie
  m$.  Then $\ad_E$ and $\ad_F$ act via the matrices $\varphi_E$ and
  $\varphi_F$ of~\eqref{eq:Sk}, so $\ad_E\ad_F$ acts as a diagonal matrix
  with entries $k$, $2(k-1)$, $3(k-2)$, \dots, $(k-1)2$, $k$ and~$0$.  As
  $\xi_A$ is in the image of~$\ad_E$, in order to have $J^2\xi_A=-\xi_A$,
  we need all the non-zero eigenvalues of~$\ad_E\ad_F$ to be~$1$.  This
  forces $k=1$.
  
  Let $\lie g(i)$ be the $i$-eigenspace of~$\ad_H$ on~$\lie g^{\mathbb C}$.
  Then $\lie p=\bigoplus_{i\geqslant 0}\lie g(i)$ is a parabolic
  subalgebra, so we may choose a Cartan subalgebra of~$\lie g^{\mathbb C}$
  lying in~$\lie p$ and a root system such that the positive root spaces
  are also in~$\lie p$.  The discussion above shows that $\ad_E$ is zero on
  all these positive root spaces, and so $E$~is a highest root vector.
  Therefore $\mathcal O=\mathcal O_{\text{min}}$.
\end{proof}

%\bibliographystyle{amsplain}
%\bibliography{papers}
\providecommand{\bysame}{\leavevmode\hbox to3em{\hrulefill}\thinspace}

\end{document}